\documentclass[a4paper,12pt]{article}

\usepackage[utf8]{inputenc}
\usepackage[english, russian]{babel}
\usepackage{amsmath, amssymb, amsthm}
\usepackage{enumerate}
 
\usepackage[colorlinks=true,linkcolor=blue,citecolor=blue]{hyperref}

\title{Об углах между линейными подпространствами в $\mathbb R^4$ и сингулярности}
\author{Чеботаренко А.О.}
\date{}

\renewcommand{\phi}{\varphi}

\renewcommand{\geq}{\geqslant}

\begin{document}

\maketitle

\hfill





Рассматривается задача о нахождении двумерных рациональных подпространств в четырёхмерном пространстве, образующих малый угол с заданным иррациональным подпространством. В этой постановке задача рассматривалась В. М. Шмидтом в 1967 году. Мы доказываем существование иррациональных подпространств, допускающих сингулярные приближения.

\section{Введение}

Феномен сингулярности в теории диофантовых приближений был обнаружен Хинчиным в его знаменитой работе
\cite{H}. 
 В частности, в этой работе он доказал, что для произвольной сколь угодно быстро убывающей к нулю функции 
$\varphi (t)$ найдутся два вещественных числа $\alpha_1$   и $\alpha_2$  такие, что   для соответствующей им функции меры иррациональности будет выполнено неравенство
$$
0< \psi_{\alpha_1,\alpha_2}(t) =\min_{x_0, x_1,x_2 \in \mathbb{Z}: 1\le \max (|x_1|,|x_2|) \le t} |x_0+x_1\alpha_1+x_2\alpha_2| \le 
\varphi(t) ,\,\,\, \forall \, t >1.
$$ 
Отметим, что согласно теореме Минковского о выпуклом теле для функции меры иррациональности для любых двух чисел 
$\alpha_1$, $\alpha_2$  всегда будет выполнено неравенство
$$
\psi_{\alpha_1,\alpha_2} (t)< t^{-2},
$$
a пара чисел  $\alpha_1$, $\alpha_2$   называется сингулярной, если 
 для любого $\varepsilon >0  $ найдется $ t_0$ такое, что для всех $ t> t_0$ выполняется неравенство 
 $\psi_{\alpha_1,\alpha_2} (t) < \varepsilon t^{-2}$.
Общий результат для произвольных систем линейных форм был задокументировал Ярником  в \cite{J}.
В настоящее время имеется много работ, посвященных вопросам существованию сингулярных векторов и 
обобщающих их объектов, а также  приложениям
(см. например, классическую книгу Касселса \cite{C}, относительно недавние работы
\cite{M}, \cite{K} и замечательный обзор \cite{G}).

В настоящей заметке мы обобщаем феномен сингулярности  для задачи, когда для заданного иррационального  линейного подпространства мы ищем рациональные подпространства, образующие с данным наименьший угол. Такая постановка задачи появилась в работе Шмидта \cite{S}. Здесь мы рассматриваем лишь простейший нетривиальный случай, в котором речь идет о приближении двумерного иррационального подпространства в $\mathbb{R}^4$ двумерным рациональным подпространством в $\mathbb{R}^4$.

\section{Подпространства и углы}

Через $||\cdot ||$  мы будем обозначать евклидову норму в $\mathbb{R}^4$, а через
 $\psi (X, Y) = \frac{\| X \land Y \|}{\|X\| \cdot \|Y\|}$ угол  (точнее синус угла) между векторами   $X, Y \in \mathbb{R}^4$. 

 Двумерное подпространство $B$ в $\mathbb{R}^4$ назовем рациональным, если у него есть базис из векторов с рациональными компонентами, множество таких подпространств обозначим $\mathfrak R_4(2)$. Высотой $H(B)$ подпространства $B$ назовём объём фундаментальной области содержащейся в нём двумерной  решётки $B\cap\mathbb{Z}^4$.

 Следуя Шмидту \cite{S}, определим первый и второй углы (точнее говоря, их синусы) между двумерными подпространствами $A$ и $B$.
  Назовём {\it первым углом} между  двумерными подпространствами $A$ и $B$ величину
 \begin {equation}\label{1}
 \psi_{1} (A, B) = \underset{X \in A\setminus{\{0\}}, Y \in B\setminus{\{0\}}}{\min} \psi (X, Y).
 \end {equation}
{\it Второй угол} $ \psi_{2} (A, B)$  между двумерными подпространствами определяется следующим образом.
Пусть минимум из определения первого угла (\ref{1})  достигается на векторах $X_{1} \in A\setminus{\{0\}}$ и $Y_{1} \in B\setminus{\{0\}}$ (если такие вектора не единственны, берём произвольные из них). Возьмём  вектор $X_{2} \in A\setminus{\{0\}},$ такой что  $X_{2} \perp X_{1}$ и вектор
$Y_{2} \in B\setminus{\{0\}}$ такой что $Y_{2} \perp Y_{1}$. Далее положим
 \begin {equation}\label{22}
 \psi_{2} (A, B) = \psi (X_{2}, Y_{2}).
 \end {equation}
 Понятно, что это определение корректно.
 \vskip+0.3cm
 
 Для двумерного подпространства $A$
 определим функцию меры иррациональности $\psi_A(t), t \in [1, +\infty)$   посредством равенства
$$
\psi_A(t) = \min\limits_{B \in \mathfrak R_4(2), H(B) \leqslant t}\psi_1(A, B)
$$
(здесь минимум берется по всем рациональным подпространствам $B$ высоты $\le t$).
Рассматривая несколько другую функцию меры иррациональности,
Шмидт \cite{S} доказал, что с некоторой абсолютной постоянной $C>0$ выполняется неравенство
$$
\min\limits_{B \in \mathfrak R_4(2), H(B) \leqslant t} (H(B))\cdot \psi_1(A, B)
< C\cdot t^{-2}.
$$
В частности, из приведенного выше неравенства следует, что для любого иррационального  двумерного подпространства $A$
 найдется бесконечно много рациональных подпространств $B$,
 таких что
 \begin{equation}\label{tree}
 \psi_1(A, B) \le C\cdot H(B)^{-3},
 \end{equation}
 или, что то же самое
 $$
 \liminf_{t\to \infty}\,\,\, t^{3}\cdot\psi_A(t) <\infty.
 $$
Точность показателя 3 в неравенстве (\ref{tree})  была недавно доказана Жозефом в прорывной работе  \cite{J}. Он показал существование
двумерных подпространств $A\subset\mathbb{R}^4$, для которых выполнено
$$
\inf_{B \in \mathfrak R_4(2)}   (H(B))^3 \cdot \psi_1(A, B) >0.
$$
Метрический результат имеется в \cite{M1}. Отметим также недавнюю работу \cite{dS} в которой, в частности, имеются результаты о втором угле ({в общем случае, об} угле с максимальным индексом).

  \vskip+0.3cm

 Пусть даны $k$ векторов $X_1, ...,X_k \in \mathbb{R}^4$. Определим  вещественную
$ k\times 4$ матрицу $M$  у которой  $j$-й столбец 
состоит из координат  вектора $X_j, j =1,2$. Для векторов $X_1,.., X_k$ рассмотрим обобщенный определитель
$$
D(X_1,...,X_k) = \sqrt{{\rm det}\,(M^t M)}.
$$
Для двумерных подпространств $A$ и $B$ в $\mathbb{R}^4$ рассмотрим произведение
$$
\Psi(A, B) = \psi_1(A,B)\psi_2(A,B),
$$
где  $\psi_1(A,B)$ и $\psi_2(A,B)$
определены в  (\ref{1}) и (\ref{22})  соответственно.
Если $X_1, X_2$ - базис в $A$, а $Y_1, Y_2$ - в $B$, то верна следующая формула (см.\cite{S}):
\begin{equation}\label{4}
\Psi(A, B) = \frac{D(X_1,X_2,Y_1,Y_2)}{D(X_1,X_2)D(Y_1,Y_2)}.
\end {equation}

Линейному двумерному  подпространству $A$ в $\mathbb {R}^4$ мы поставим в соответствие однородные  плюккеровы коодинаты
\begin{equation}\label{oo}
(\eta_1:\eta_2:\eta_3:\eta_4:\eta_5:\eta_6),
\end{equation}
определяемые как набор всех $2\times 2$ миноров $2\times 4$ матрицы, составленной из координат векторов какого-либо базиса в $A$ 
 (по поводу определения и свойств плюккеровых координат см. 
 \cite{hp}; там они называются "грассмановы координаты"). 
Вообще говоря, компоненты $\eta_j, 1\le j \le 6$
определены с точностью до пропорциональности, но иногда нам будут
нужны конкретные представитель однородного набора (\ref{oo}).
В этом случае мы будем  писать их через запятую:
$$
(\eta_1,\eta_2,\eta_3,\eta_4\eta_5,\eta_6).
$$ 
Множество всех двумерных линейных подпространств в $\mathbb {R}^4$ 
имеет структуру многообразия, которое назвыется грассманианом. 
Грассманиан вкладывается в $\mathbb {PR}^5$ с помощью плюккеровых координат,  для которых выполнено следующее соотношение
\begin {equation}\label{gg}
x_1x_6 - x_2x_5 + x_3x_4 = 0.
\end {equation}
Поверхность в проективном пространстве  $\mathbb {PR}^5$, точки которой удовлетворяют 
соотношению (\ref{gg}) мы будем обозначать $\mathfrak G$ и будем отождествлять ее с грассманианом, то есть, со множеством всех двумерных линейных подпространств в  $\mathbb {R}^4$ .
На   $\mathfrak G$  имеется  топология индуцированная   стандартной топологией в $\mathbb {PR}^5$. Теперь, когда речь будет идти об открытых и замкнутых подмножествах $\mathfrak G$, по умолчанию будем считать, что у нас имеется в виду именно эта топология.

 Пусть теперь в $A$ и $B$ взяты ортонормированные базисы, на векторах которых определяются первый и второй углы между $A$ и $B$.
 Рассмотрим плюккеровы координаты
  $(\eta_1, \eta_2, \eta_3, \eta_4, \eta_5, \eta_6)$ подпространства  $A$ в базисе $(X_1,X_2)$
  и  плюккеровы координаты $(\xi_1, \xi_2, \xi_3, \xi_4, \xi_5, \xi_6)$ подпространства  $B$ в базисе $(Y_1,Y_2)$.  
 Тогда
$$
D(X_1,X_2) = \sqrt{\eta_1^2 + \eta_2^2 + \eta_3^2 + \eta_4^2 + \eta_5^2 + \eta_6^2} = 1 ,
$$
$$
D(Y_1,Y_2) = \sqrt{\xi_1^2 + \xi_2^2 + \xi_3^2 + \xi_4^2 + \xi_5^2 + \xi_6^2} = 1 ,
$$
$$
D(X_1,X_2,Y_1,Y_2) = |\xi_1\eta_6 + \xi_6\eta_1 - \xi_2\eta_5 - \xi_5\eta_2 + \xi_3\eta_4 + \xi_4\eta_3|.
$$
Таким образом из $(\ref{4})$ получаем
\begin {equation}\label{7}
\Psi(A, B) = |\xi_1\eta_6 + \xi_6\eta_1 - \xi_2\eta_5 - \xi_5\eta_2 + \xi_3\eta_4 + \xi_4\eta_3|.
\end {equation}

\section{Сингулярность}
Теперь мы сформулируем основной результат настоящей работы.

 \vskip+0.3cm
{\textbf{Теорема 1.}}{\it
Пусть $\phi(t)$ - произвольная строго положительная убывающая функция на $[1,+\infty)$. Тогда существует такое подпространство  $\mathcal A$ и $t_0 > 1$, что для $\forall t > t_0$ выполнено
$$
0 < \psi_\mathcal A(t) < \phi(t).
$$
}
\vskip+0.3cm

Опишем схему доказательства теоремы 1.

Будем индуктивно строить последовательность рациональных подпространств $Q_i$ с высотами $H_i$, такими что
$$
1) \,\,H_i < H_{i+1}, \forall i \in \mathbb N;\,\,\,\,\,\,\,\,\,\,\,
2) \,\,\lim \limits_{n \to \infty} H_i = \infty .
$$
С помощью этих подпространств мы построим последовательность замкнутых вложенных множеств 
\begin{equation}\label{re}
 \mathfrak G \supset \mathfrak A_1 \supset \mathfrak A_2 \supset \mathfrak A_3 ...   
\end{equation}
таких, что для $\forall A \in \mathfrak A_i, \forall i \geq 2$ выполнено 
\begin {equation}\label{ro}
1) \,\,\psi_1(A, Q_{i-1}) < \phi(H_{i}) ;\,\,\,\,\,\,\,\,\,\,\,
2)\,\, \psi_A(H_{i-1}) > 0.
\end {equation}
Так как грассманиан    $\mathfrak G$  компактен,  выполняется соотношение
$$
\bigcap_{n=1}^\infty \mathfrak A_n \neq  \varnothing.
$$
Рассмотрим произвольное подпространство $\mathcal A$, принадлежащее этому пересечению. Покажем, что для него выполнено условие теоремы. Действительно, для него выполнено
$$
0 < \psi_\mathcal A(H_i) < \phi(H_{i+1}), \forall i \in \mathbb N.
$$
Легко видеть, что функция $\psi_\mathcal A(t)$ невозрастающая, поэтому, так как $\lim \limits_{n \to \infty} H_i = \infty$, то $\psi_\mathcal A(t) > 0$ для любого $ t \in [1, \infty)$. Пусть $t \in [H_i, H_i+1)$, тогда, в силу убывания $\phi (t)$ имеем
$$
\psi_\mathcal A(t) \leqslant \psi_\mathcal A(H_i) < \phi(H_{i+1}) < \phi(t).
$$
Взяв теперь $t_0 = H_1$, получаем требуемое.

\vskip+0.3cm
Последовательность множеств (\ref{re}), подпространства в которых удовлетворяют ({\ref{ro}), мы построим в пункте 5,
а в следующем пункте 4 мы сформулируем и докажем простейшие вспомогательные утверждения.

\section{Вспомогательные утверждения}\label{sec:intro}

Для подпространства $A$ рассмотрим множество $\mathfrak V(A)$ подпространств $B$, таких что $\Psi (A,B) = 0$. Пусть $(\eta_1: \eta_2: \eta_3: \eta_4: \eta_5: \eta_6)$ - плюккеровы координаты $A$, а $(\xi_1: \xi_2: \xi_3: \xi_4: \xi_5: \xi_6)$ - плюккеровы координаты $B$. Тогда $B \in \mathfrak V(A)$ тогда и только тогда, когда 
$$
\xi_1\eta_6 + \xi_6\eta_1 - \xi_2\eta_5 - \xi_5\eta_2 + \xi_3\eta_4 + \xi_4\eta_3 = 0.
$$

Здесь стоит сделать замечание, что имеет место следующая цепочка эквивалентностей
$$
B \in \mathfrak V(A) \Leftrightarrow A \in \mathfrak V(B) \Leftrightarrow A \cap B \neq \varnothing.
$$

В этом пункте мы докажем два совсем простых вспомогательных утверждения.
Для ненулевого вектора $x \in \mathbb R^4$ рассмотрим порожденное им одномерное линейное подпространство 
$\ell (x)$  и ведём отображение $g_x: \mathbb R^4 \setminus \ell (x) \rightarrow \mathfrak G$, которое сопоставляет вектору $y \in \mathbb R^4$ точку в $\mathfrak G$, соответствующую двумерному подпространству, натянотому на вектора $x$ и $y$. Легко видеть, что для любого $x \neq 0$
отображение  $g_x$ является непрерывным.

\vskip+0.3cm

{\textbf{Лемма 1. }}\, {\it
Пусть дано подпространство $Q \in \mathfrak G$, а также дано конечное число подпространств $P_1, P_2, ..., P_n \in \mathfrak G$ и $P_i \neq Q$ для любого $i$ . Тогда для любого открытого множества $\mathfrak U$
c условием
$ \mathfrak U \cap \mathfrak V(A) \neq \varnothing$ найдётся открытое множество $\mathfrak U_0 \subset \mathfrak U$, такое что $\mathfrak U_0 \cap \mathfrak V(A) \neq \varnothing$ и  $\mathfrak U_0 \cap \mathfrak V(P_i) = \varnothing$ для любого $i = 1,2,...,n$.}
\vskip+0.3cm

Доказательство.

Ясно, что достаточно доказать лемму при $n=1$.
Построим открытое множество $\mathfrak U_0 \subset \mathfrak U$, такое что $\mathfrak U_0 \cap \mathfrak V(Q) \neq \varnothing$, но $\mathfrak U_0 \cap \mathfrak V(P_1) = \varnothing$. Возьмём произвольное $A \in \mathfrak U \cap \mathfrak V(Q)$, такое что $A \neq P_1$   и $A \neq Q$. 

Разберем два случая.

Случай 1$^0$.
Если $A \not\in \mathfrak V(P_1)$, то в силу того что $\mathfrak G \setminus \mathfrak V(P_1)$ открыто, найдётся открытая окрестность   $A$, удовлетворяющая всем требуемым  условиям, её и возьмём в качестве $\mathfrak U_0$.

Случай 2$^0$.
Если же 
 $A\in \mathfrak V(P_1)$, 
 то мы возьмем  ненулевой вектор $e \in Q \cap A$
 (пересечение  $ Q \cap A$ непусто, так как
 $\psi_1 (A, Q) = 0$).
Если $e \in P_1$ 
(в этом случае $ P_1\cap A = \ell(e)$ является одномерным подпрпостранством порожденным вектором $e$)
то мы берем произвольный произвольный вектор $f \in  A$ не пропорциональный вектору $e$.
В силу непрерывности функции $g_f$, мы можем найти вектор $e' \in Q$, такой что $A' = g_f (e') \in \mathfrak U$. Но $e' \in A' \cap Q$, поэтому  $A' \in  \mathfrak U \cap \mathfrak V(Q)$. Теперь для подпространства $A'$ общий вектор с подпространством $Q$ уже не будет принадлежать подпространству $P_1$, и вместо $A$ мы можем продолжить работать с $A'$ (штрих при $A$ мы в дальнейшем писать не будем). Итак, можно считать, что $e\in Q\cap A$  и $e \not\in  P_1$. Тогда все двумерные подпространства, содержащие $e$ и пересекающие $P_1$ лежат в трёхмерном подпространстве 
порожденном двумерным подпространством $P_1 $ и вектором $e$.
Это подпространство мы обозначим через
$\langle P_1, e\rangle$.

Отметим, что поскольку $A$  порождено векторами $e$ и $f$, выполняется равенство
$g_e(f) = A$.
 В силу непрерывности функции $g_e$ и нигде не плотности  подпространства  $\langle P_1, e\rangle$ в объемлющем пространстве  $\mathbb R^4$, найдётся достаточно близкий к $f$ вектор $f''$, не лежащий в этом трёхмерном подпространстве, такой что $A'' = g_e(f'') \in \mathfrak U$. При этом, так как $e \in A''$, то $A'' \in \mathfrak V(Q)$.
 Так как $f''\not\in \langle P_1, e\rangle $, то, тем более, $A'' \not\subset \langle P_1, e\rangle $, и  мы заключаем, что $A'' \not\in \mathfrak V(P_1)$.
  Теперь мы можем действовать так же, как в случае 1$^0$, заменив $A$ на $A''$.$\Box$

 \vskip+0.3cm

{\textbf{Лемма 2. }}\, {\it
Обозначим через $\mathfrak Q_A$ множество рациональных подпространств, лежащих в $ \mathfrak V(A)$. Рассмотрим в $\mathfrak V(A)$ топологию, индуцированную топологией в  $\mathfrak G$. Тогда, если $A \in \mathfrak R_4(2)$, то $\mathfrak Q_A$ всюду плотно в этой топологии.}

\vskip+0.3cm

Доказательство.

Достаточно показать, что для любого открытого множества $\mathfrak U \subset \mathfrak G$, такого что $\mathfrak U \cap \mathfrak V(A) \neq \varnothing$, найдётся рациональное подпространство $Q \in U \cap \mathfrak V(A)$. Пусть $B$ - произвольное подпространство, лежащее в $\mathfrak U \cap \mathfrak V(A)$, тогда существует ненулевой вектор $e \in A \cap B$. Рассмотрим ненулевой вектор  $f \in B$, такой что  $f \neq e$. Так как рациональные векторы всюду плотны в подпространстве $A$, то найдётся рациональный вектор $q_1 \in A$, такой что $B' = g_f(q_1) \in \mathfrak U$. Пусть вектор $h \in B'$ и $h \neq q_1$. Тогда, так как $\mathfrak U$ открыто и $B' \in \mathfrak U$, а множество рациональных векторов всюду плотно в $\mathbb R^4$, найдётся достаточно близкий к $h$ рациональный вектор $q_2$, такой что $Q = g_{q_1}(q_2) \in \mathfrak U$. Тогда $Q \in \mathfrak R_4(2)$, и при этом $Q \in \mathfrak V(A)$, так как $q_1$ - общий вектор $Q$ и $A$. Отсюда следует, что $\mathfrak Q_A$ всюду плотно в $\mathfrak V(A)$.$\Box$

\section{Построение последовательности рациональных подпространств и множеств $\mathfrak A_i$}

 Введём функцию $ \Psi_A(B): \mathfrak G \rightarrow [0, 1]$ следующим образом:
$$
\Psi_A (B): B \mapsto \Psi (A,B).
$$
 Легко видеть {{из (\ref{7})}}, что  для любого $A$ из $\mathfrak G$ функция $\Psi_A(B)$ непрерывна. Следовательно, $\mathfrak V(A)=\Psi_A^{-1} (\{ 0\})$ будет замкнуто для любого $A$.
Будем обозначать через  ${\rm int}\,(\mathfrak A)$ - множество внутренних точек множества $\mathfrak A$.
Как уже говорилось, множества $\mathfrak{A}_i$, удовлетворяющие (\ref{re}) и (\ref{ro}), мы будем строить индуктивным образом, причем мы потребуем, чтобы еще дополнительно выполнялось условие
\begin{equation}\label{int}
\mathfrak V(Q_i) \cap {\rm int}\,(\mathfrak A_i) \neq \varnothing.
\end{equation}
Положим $\mathfrak A_1 = \mathfrak G$. А в качестве $Q_1$ возьмём произвольное рациональное подпространство. Пусть теперь построены удовлетворяющие нашим условиям множества $\mathfrak A_1, \mathfrak A_2, \mathfrak A_3, ..., \mathfrak A_i$ и соответствующие  им рациональные подпространства $Q_1, Q_2, Q_3, ..., Q_i$. В силу леммы $2$, и того что $\mathfrak V(Q_i)$ пересекает внутренность $\mathfrak A_i$, найдётся рациональное подпространство
$$
Q \in \mathfrak V(A_i) \bigcap {\rm int}\,(\mathfrak A_i) .
$$
Понятно, что мы можем взять $Q$ отличным от $Q_i$, тогда $\psi_2(Q_{i}, Q) > 0$. Снова  воспользовавшись леммой $2$, мы можем взять рациональное подпространство $Q_{i+1} \in \mathfrak V(Q)$, такое что $H_{i+1} := H(Q_{i+1}) > H_i$. Рассмотрим множество  $\mathfrak F (Q_i)$ подпространств $A$, таких что
\begin{equation}\label{u1}
\psi_2 (A, Q_i) > \frac{\psi_2(Q_{i}, Q)}{2} .
\end{equation}
Оно, очевидно, открыто. Также открыто будет множество 
$\mathfrak V_{i}$ таких подпространств $A$, что
\begin{equation}\label{u2}
\Psi (A, Q_{i}) < \frac{\psi_2(Q_{i}, Q)}{2} \phi(H_{i+1}).
\end{equation}
Заметим, что тогда множество $\mathfrak U(Q)=\mathfrak V_{i} \cap {\rm int} (\mathfrak A_i) \cap  \mathfrak F (Q_i)$ будет открыто, и, при этом, не пусто, так как в нём лежит $Q$. Согласно 
(\ref{u1}) и  (\ref{u2}) для любого $A \in \mathfrak U(Q)$ будет выполнено
\begin{equation}\label{17}
\psi_1 (A, Q_{i}) = \frac{\Psi (A, Q_i)}{\psi_2(A,Q_{i})} < \frac{\frac{\psi_2(Q,Q_{i})}{2}\phi(H_{i+1})}{\frac{\psi_2(Q,Q_{i})}{2}} = \phi(H_{i+1}).
\end{equation}
Рассмотрим множество рациональных подпространств
$$
\mathfrak R_{i+1} :=  \{R \in \mathfrak R_4(2), H_{i-1} < H(R) \leqslant H_{i}\}
$$
(полагаем $H_{i-1} = 0$ при $i=1$).
В силу того что рациональных подпространств с ограниченной высотой конечное число, множество $\mathfrak R_{i+1}$ конечно, а значит мы можем применить лемму $1$ для подпространства $Q_{i+1}$, набора подпространств из $\mathfrak R_{i+1}$ и открытого множества $\mathfrak U(Q)$. Согласно  этой лемме найдётся такое открытое множество $\mathfrak U_0(Q)$, что 
$$
1)\,\, \mathfrak U_0(Q) \cap \mathfrak V(Q_i) \neq \varnothing;\,\,\,\,\,\,\,\,\,\,
2)\,\, \mathfrak U_0(Q) \cap \mathfrak V(R) = \varnothing ,\,\,\,\forall R \in \mathfrak R_{i+1}.
$$
Понятно, что в  $\mathfrak U_0(Q)$ найдётся замкнутое множество, такое что множество его внутренних точек всё ещё пересекает $\mathfrak V(Q_{i+1})$. Возмём это множество в качестве $\mathfrak A_{i+1}$, покажем, что для него выполняются все требуемые условия. Во-первых,
\begin{equation}\label{01}
{\rm int}\,(\mathfrak A_{i+1}) \cap \mathfrak V(Q_i) \neq \varnothing.
\end{equation}
Во-вторых, в силу того что $\mathfrak A_{i+1} \subset \mathfrak A_{i}$, из требований на $\mathfrak A_{i}$ видим, что выполнено
$$
\psi_A (H_{i-1}) > 0,\,\,\, \forall A \in  \mathfrak A_{i+1}.
$$
Из построения следует, что любое подпространство из $\mathfrak A_{i+1}$ не пересекает никакое рациональное подпространство с высотой $H$, такой что $H_{i-1} < H \leqslant H_i$. Значит, на самом деле выполняется
\begin{equation}\label{02}
\psi_A (H_{i}) > 0,\,\,\, \forall A \in  \mathfrak A_{i+1}.
\end{equation}
И, наконец, так как $\mathfrak A_{i+1} \subset \mathfrak U(Q)$, то из (\ref{17}) получаем
\begin{equation}\label{03}
\psi_1 (A, Q_{i}) < \phi(H_{i+1}), \,\,\, \forall A \in \mathfrak A_{i+1}.
\end{equation}
Итак, соотношения (\ref{01}), (\ref{02}) и (\ref{03}) обеспечивают условия (\ref{ro}) и (\ref{int}) для следующего шага индукции.
Это завершает доказательство.$\Box$

\vskip+0.5cm
Работа выполнена при финансовой поддержке фонда Базис, грант № 21-7-1-33-1.

\vskip+0.5cm

\end{document}